\documentclass[11pt,a4paper]{article}
\usepackage{a4}
\usepackage{amssymb,latexsym}
\def\eqnreset{\setcounter{equation}{0}}
\def\eqsection#1{\section{#1}\eqnreset}

\newtheorem{Thm}{Theorem}[section]
\newtheorem{Defi}[Thm]{Definition}
\newtheorem{Cor}[Thm]{Corollary}
\newtheorem{Lemma}[Thm]{Lemma}
\newtheorem{Prop}[Thm]{Proposition}
\newtheorem{Rem}[Thm]{Remark}
\newtheorem{Conj}[Thm]{Conjecture}
\newtheorem{Prelim}[Thm]{Preliminary}

\newenvironment{thm}[0]{\begin{Thm}\noindent}%
{\end{Thm}}
\newenvironment{defi}[0]{\begin{Defi}\noindent\rm}%
{\end{Defi}}
{\end{Cor}}
\newenvironment{lemma}[0]{\begin{Lemma}\noindent}%
{\end{Lemma}}
{\end{Prop}}
\newenvironment{rem}[0]{\begin{Rem}\noindent\rm}%
{\end{Rem}}
{\end{Conj}}
{\end{Prelim}}

\def\qed{~\hfill$\square$\medbreak}

\def\text#1{\;\;\;\;{\rm \hbox{#1}}\;\;\;\;}
\def\qquad{\quad\quad}

\def\msy#1{{\mathbb #1}}
\def\C{{\msy C}}
\def\N{{\msy N}}
\def\Z{{\msy Z}}
\def\real{{\msy R}}
\def\D{{\msy D}}
\def\T{{\msy T}}

\def\ga{\alpha}

\def\gl{\lambda}

\def\gL{\Lambda}
\def\gS{\Sigma}

\def\frak#1{\mathfrak #1}
\def\fa{{\frak a}}

\def\fg{{\frak g}}
\def\fh{{\frak h}}

\def\fk{{\frak k}}

\def\fq{{\frak q}}

\def\fu{{\frak u}}

\def\to{\rightarrow}
\def\Re{{\rm Re}\,}
\def\Im{{\rm Im}\,}

\def\Ad{{\rm Ad}}

\def\cU{{\mathcal U}}

\def\rmO{{\rm O}}
\def\rmS{{\rm S}}
\def\rmU{{\rm U}}
\def\SU{\rmS\rmU}
\def\SO{\rmS\rmO}
\def\PW{{\rm PW}}
\def\Exp{{\rm Exp}}
\begin{document}
\setcounter{section}{0}
\title{A local Paley--Wiener theorem for\\ compact symmetric spaces}
\author{Gestur \'Olafsson%
\footnote{Research supported by NSF grant DMS-0402068}
\ and Henrik Schlichtkrull}
\date{26 April 2007}
\maketitle
\begin{abstract}
The Fourier coefficients of a smooth $K$-invariant function on 
a compact symmetric space $M=U/K$ 
are given by integration of the function against the spherical functions.
For functions with support in a neighborhood of the origin, 
we describe the size of the support by means of the exponential type
of a holomorphic extension of the Fourier coefficients.
\footnote{2000 Mathematics Subject Classification: 33C55, 43A85, 53C35}
\end{abstract}
\noindent
\eqsection{Introduction}

The classical Paley-Wiener theorem (also called the
Paley-Wiener-Schwartz theorem)
describes the image by the Fourier
transform of the space of compactly supported smooth
functions on $\real^n$. The theorem was generalized to
Riemannian symmetric spaces 
of the noncompact type by Helgason and Gangolli 
(see \cite{GGA}, Thm.\ IV,7.1, \cite{Gang}), 
to semisimple Lie groups by Arthur (see \cite{Arthur}),
and to pseudo-Riemannian
reductive symmetric spaces by van den Ban and Schlichtkrull
(see \cite{BanS}). More precisely, these theorems describe
the Fourier image of the space of functions supported in a 
(generalized) ball of a given size. The image space
consists of holomorphic (in the pseudo-Riemannian case,
meromorphic) functions with exponential growth,
and the size of the ball is reflected in the exponent
of the exponential growth estimate.

In this paper we present an analogue of these theorems
for Riemannian symmetric spaces of the compact type. 
Obviously the compact support is trivial in this case,
and the important issue is the determination of 
the {\it size} of the support of a smooth function from the
growth property of its Fourier transform. Let us illustrate
this by recalling the corresponding result for Fourier 
series. Consider a smooth $2\pi$-periodic function 
$f\colon \T=\real/2\pi\Z\to\C$, and suppose that
$f$ has support in $[-r,r]+2\pi\Z$, where $0<r<\pi$.
We denote the space of such functions by $C^\infty_r(\T)$.
The Fourier transform of $f$ is the Fourier coefficient
map $n\mapsto \hat f(n)$ on $\Z$, where
$$\hat f(n)=\frac1{2\pi}
\int_{-\pi}^\pi f(e^{it}) e^{-int}\,dt,$$
and it extends to a holomorphic function on $\C$,
defined by the same formula with $n$ replaced by
$\gl\in\C$. By the classical Paley-Wiener theorem 
for $\real$ this holomorphic extension has at most
exponential growth of type $r$, and every holomorphic 
function on $\C$ of this type arises in this fashion
from a unique function $f\in C^\infty_r(\T)$.
It is this 'local' Paley-Wiener theorem for $\T$ that we 
generalize to an arbitrary Riemannian symmetric space
$M$ of the compact type. We consider spherical functions
on $M$, and the relevant transform is the
spherical Fourier transform.

The theorem presented here is known in some particular cases.
In particular, it is known in the case of a compact
Lie group $U$, viewed as a symmetric space for the 
product group $U\times U$ with the left$\times$right action.
In this case, the theorem was obtained by Gonzalez (see 
\cite{Gon}) by a simple reduction to the Euclidean case
by means of the Weyl character formula. This result of
Gonzalez plays a crucial role in our proof, and it is
recalled in Section \ref{s: central} below.
Other cases in which the theorem is known, are as follows.

If the symmetric space has rank one, the spherical
Fourier transform can be expressed in terms of a
Jacobi transform, for which the Paley-Wiener
theorem has been obtained by Koornwinder (see \cite{TomK} p.\ 158). 
As an example, we treat the special case $S^2=\SU(2)/\SO(2)$ 
in the final section of this paper. In this case, the
theorem of Koornwinder is due to Beurling (unpublished,
see \cite{TomK}). 

If the symmetric space is of even multiplicity type,
the local Paley-Wiener theorem has been achieved by Branson, 
\'Olafsson and Pasquale (see \cite{BOP}) by application
of a holomorphic version of Opdam's differential shift operators 
(developed in \cite{Opd}, \cite{OlPa}).
The method  is strongly dependent on the assumption that
the multiplicities are even. The theorem of Gonzalez is a 
particular case.

Finally, the theorem was obtained recently by Camporesi
for the complex Grassmann manifolds by reduction to the
rank one case, see \cite{Camp}.

We shall now give a brief outline of the paper.
In Sections \ref{s: notation} 
and \ref{s: Fourier} we introduce the basic notations.
In Section \ref{s: main thm} we define the relevant 
Paley-Wiener space and state the main theorem, that
the Fourier transform is bijective onto this space.
The proof, that it maps into the space is given in
Section \ref{s: Opdam}. Here we rely on work of Opdam \cite{Opd}.
The theorem of Gonzalez, mentioned above, is recalled
in Section \ref{s: central}, and the central argument 
of the present paper, establishing surjectivity, 
is given in the following Sections 
\ref{s: K-invariant}-\ref{s: surjective}. An important
ingredient is a result of Rais from \cite{Rais}, which
has previously been applied in similar situations
by Clozel and Delorme, \cite{CD1} Lemma 7, and by 
Flensted-Jensen
\cite{FJ} p.\ 30. Finally, in Section \ref{s: sphere} we treat
$S^2$ as an example.

The result of \cite{BOP} has been
generalized to the Jacobi transform
associated to a root system with 
a multiplicity function which is even,
but not necessarily related to a symmetric space
(see \cite{BOP2}). 
For the method of the present paper the geometry of
the symmetric space is crucial, especially in Lemma
\ref{l: support}, and we do not see how to
generalize in this direction.

\eqsection{Basic notation}
\label{s: notation}
Let $M$ be a Riemannian symmetric space of the compact type.
We can write $M$ as a homogeneous space
$M=U/K$, where $U$ is a connected compact 
semisimple Lie group which acts isometrically on $M$, and 
$K$ a closed subgroup with the property that  
$U^\theta_0\subset K\subset U^\theta$ for an
involution $\theta$ of $U$. Here $U^\theta$ denotes the 
subgroup of $\theta$-fixed points, and 
$U^\theta_0$ its identity component. It should be emphasized
that the pair $(U,K)$ is in general not uniquely
determined by $M$ (see \cite{Sig}, Ch. VII).

Let $\fu$ denote the Lie algebra of $U$. We denote
the involution of $\fu$ corresponding to $\theta$ by the 
same symbol. Let $\fu=\fk\oplus\fq$ be the corresponding 
Cartan decomposition, then $\fk$ is the Lie algebra of
$K$, and $\fq$ can be identified with the tangent space
$T_oM$, where $o=eK\in M$ is the origin.

Recalling that the Killing form $B(X,Y)$ on $\fu$ is 
negative definite, let $\langle\,\cdot\,,\,\cdot\,\rangle$
be the inner product on $\fu$ defined by
$\langle X,Y\rangle=-B(X,Y)$. 
Then $\fk$ and $\fq$ are orthogonal subspaces. 
We assume that the Riemannian metric $g$ of $M$ is
normalized such that it agrees with $\langle\,\cdot\,,\,\cdot\,\rangle$
on $\fq=T_oM$. 

We denote by $\exp$ the exponential map $\fu\to U$ 
(which is surjective), and by $\Exp$ the map $\fq\to M$ given by
$\Exp(X)=\exp(X) \cdot o$.
By identification of $\fq$ with the tangent space 
$T_oM$, we thus identify $\Exp$ with the exponential map
associated to the Riemannian connection.

The inner product on $\fu$
determines an inner product on the dual space $\fu^*$
in a canonical fashion. Furthermore, these inner products
have complex bilinear extensions to the complexifications
$\fu_\C$ and $\fu_\C^*$. All these bilinear forms are
denoted by the same symbol $\langle\,\cdot\,,\,\cdot\,\rangle$.

Let $\fa\subset\fq$ be a maximal abelian subspace,
$\fa^*$ its dual space, and $\fa^*_\C$ the complexified
dual space. Let $\Sigma$ denote the set of non-zero 
(restricted) roots of 
$\fu$ with respect to $\fa$, then $\Sigma\subset\fa^*_\C$ and all
the elements of $\Sigma$ are purely imaginary on $\fa$.
The multiplicity of a root $\alpha\in\Sigma$ is denoted $m_\alpha$.
The corresponding Weyl group, generated by the reflections
in the roots, is denoted $W$. Recall that it is
naturally isomorphic with the factor group $N_K(\fa)/Z_K(\fa)$ of the
normalizer and the centralizer of $\fa$ in $K$
(see \cite{Sig}, Cor.\ VII.2.13).

\eqsection{Fourier series}
\label{s: Fourier}

Let $(\pi,V)$ be an irreducible unitary representation of
$U$, and let 
$$V^K=\{v\in V\mid \forall k\in K: \pi(k)v=v\},$$
then $V^K$ is either $0$ or 1-dimensional. In the latter
case $\pi$ is said to be a $K$-{\it spherical} representation.

Let $\fh\subset\fu$ be a Cartan subalgebra containing $\fa$,
then $\fh=\fh_m\oplus\fa$, where $\fh_m=\fh\cap\fk$.
Let $\Delta$ denote the set of roots of $\fu$ with respect
to $\fh$, then $\Sigma$ is exactly the set of non-zero restrictions
to $\fa$ of elements of $\Delta$. We fix
a set $\Sigma^+\subset\Sigma$ of positive restricted roots,
and a compatible set $\Delta^+\subset\Delta$ of positive roots.
The set of dominant integral linear functionals on $\fh$ is
$$\gL^+(\fh)=\{\gl\in\fh_\C^*\mid \forall\alpha\in\Delta^+:
\frac{2\langle\gl,\alpha\rangle}{\langle\alpha,\alpha\rangle}
\in\Z^+\},$$
where $\Z^+=\{0,1,2,\dots\}$.
Notice that since $\fu$ is compact,
all elements of $\Delta$ and $\Lambda^+(\fh)$ 
take purely imaginary values on on $\fh$.

Let $\gL^+(U)\subset\fh^*$ denote the set
of highest weights of irreducible representations
of $U$, then $\gL^+(U)\subset\gL^+(\fh)$ with equality
if and only if $U$ is simply connected. Let
$\gL^+_K(U)$ denote the subset of $\gL^+(U)$
which corresponds to $K$-spherical representations.
We recall the following identification of 
$\gL^+_K(U)$,
due to Helgason (see \cite{GGA}, p.\ 535).

\begin{thm} 
\label{t: Helgason classification}
Let $\gl\in\gL^+(U)$.
Then $\gl\in\gL^+_K(U)$ if and only if
$\gl|_{\fh_m}=0$ and the restriction $\mu=\gl|_{\fa}$
satisfies
\begin{equation}
\label{e: Helgason condition}
\frac{\langle\mu,\alpha\rangle}{\langle\alpha,\alpha\rangle}
\in\Z^+,
\end{equation}
for all $\alpha\in\Sigma^+$.

Furthermore, if $\mu\in\fa^*$ satisfies
{\rm (\ref{e: Helgason condition})}
for all $\alpha\in\Sigma^+$, then the element
$\gl\in\fh^*_\C$ defined by
$\gl|_{\fh_m}=0$ and $\gl|_{\fa}=\mu$
belongs to $\gL^+(\fh)$. If this element
$\gl$ belongs to $\gL^+(U)$, then it belongs to $\gL^+_K(U)$.
\end{thm}

Let $\gL^+(U/K)$ denote the set of restrictions $\mu=\gl|_{\fa}$
where $\gl\in\gL^+_K(U)$, according to the preceding theorem
this set is in bijective correspondence with $\gL^+_K(U)$.
For each $\mu\in\gL^+(U/K)$ we fix an irreducible unitary
representation
$(\pi_\mu,V_\mu)$ of $U$ with highest weight $\gl$, and we fix
a unit vector $e_\mu\in V_\mu^K$.
The  {\it spherical function} on $U/K$ associated
with $\mu$ is the matrix coefficient
$$\psi_\mu(x)=(\pi_\mu(x)e_\mu,e_\mu), \quad x\in U,$$
viewed as a function on $U/K$.
It is $K$-invariant on both sides, and it is
is independent of the choice of the unit vector $e_\mu$. The
{\it spherical Fourier transform} of a continuous 
$K$-invariant function $f$ on
$M=U/K$ is the function $\tilde f$ on $\gL^+(U/K)$ defined by
$$\tilde f(\mu)=\int_M f(x)\overline{\psi_\mu(x)}\,dx,$$
where $dx$ is the Riemannian measure on $M$, normalized
with total measure 1. Notice that $\overline{\psi_\mu(gK)}
=\psi_\mu(g^{-1}K)$ for $g\in U$, since $\pi_\mu$ is unitary.
The {\it spherical Fourier series}
for $f$ is the series given by
\begin{equation}
\label{e: U/K Fourier series}
\sum_{\mu\in\gL^+(U/K)} d(\mu)\tilde f(\mu)\psi_\mu
\end{equation}
where $d(\mu)=\dim V_\mu$. The Fourier series converges to $f$ 
in $L^2$ and, if $f$ is smooth,
absolutely and uniformly (see \cite{GGA}, p.\ 538).

Furthermore, $f$ is smooth if and only if the Fourier 
transform $\tilde f$ is {\it rapidly decreasing}, that is,
for each $k\in\N$ there exists a constant $C_k$ such
that
$$|\tilde f(\mu)|\leq C(1+\|\mu\|)^{-k}$$
for all $\mu\in\gL^+(U/K)$ (see \cite{Sugiura}).

\eqsection{Main theorem}
\label{s: main thm}

For each $r>0$ we denote by $B_r(0)$ the open ball in $\fq$ 
centered at $0$ and with radius $r$. The exponential image
$\Exp B_r(0)$ is the ball in $M$, centered at the origin 
and of radius $r$. 
Let $\bar B_r(0)$ and $\Exp \bar B_r(0)$
denote the corresponding closed balls.
We denote by $C^\infty_r(U/K)^K$
the space of $K$-invariant smooth functions on $M=U/K$
supported in $\Exp \bar B_r(0)$.

Let $\rho=\frac12\sum_{\alpha\in\Sigma^+} m_\alpha\alpha\in\fa^*_\C$.

\begin{defi} 
\label{d: PW space}
(Paley-Wiener space)
For $r>0$ let
$\PW_r(\fa)$
denote the space of holomorphic functions $\varphi$ on $\fa_\C^*$
satisfying the following.

\item{(a)} For each $k\in\N$ there exists a constant $C_k>0$ such that
$$|\varphi(\gl)|\leq C_k(1+\|\gl\|)^{-k} e^{r\|\Re\gl\|}$$
for all $\gl\in\fa_\C^*.$

\item{(b)} $\varphi(w(\gl+\rho)-\rho)=\varphi(\gl)$ for
all $w\in W$, $\gl\in\fa_\C^*$.

\end{defi}

We can now state the main theorem.

\begin{thm}
{\rm (The local Paley--Wiener theorem)}
\label{t: PW}
There exists $R>0$ such that the following holds
for each $0<r<R$.
\smallskip

{\rm (i)} Let $f\in C^\infty_r(U/K)^M$.
Then the Fourier transform $\tilde f\colon \Lambda^+(U/K)\to\C$
of~$f$ extends to a function in $\PW_r(\fa)$.

{\rm (ii)} Let $\varphi\in\PW_r(\fa)$. There exists
a unique function $f\in C^\infty_r(U/K)^K$ such that
$\tilde f(\mu)=\varphi(\mu)$ for all $\mu\in \Lambda^+(U/K)$.

{\rm (iii)} The functions in the Paley-Wiener space 
$\PW_r(\fa)$ are uniquely determined
by their values on $\Lambda^+(U/K)$.

\smallskip
Thus the Fourier transform followed by the extension gives a bijection
$$C^\infty_r(U/K)^K\to\PW_r(\fa).$$
\end{thm}

\begin{rem}\label{r: best R} 
It would be reasonable to expect the theorem above
to hold with $R$ equal to the injectivity radius of $M$, 
that is, the supremum of the values $r$ for which the
restriction of $\Exp$ to $B_r(0)$ is a
diffeomorphism onto its image. We have not been able to 
establish that. It should be noted that parts (i)-(iii)
of the above theorem may be valid with different values of
$R$. In fact, it can be seen from the proofs below, that if 
$U$ is simply connected, then part (i) of the theorem is valid 
with $R$ equal to half the injectivity radius of $M$. Furthermore,
part (ii) will be established with $R$ equal to the injectivity 
radius of $U$. For part (iii) we need a possibly smaller value of $R$.
\end{rem}

\eqsection{The invariant differential operators}
\label{s: D(U/K)}
Let $\D(U/K)$ denote the algebra of $U$-invariant differential
operators on $U/K$. It is commutative (see \cite{GGA}, Cor.\ II.5.4).
Recall that the {\it Harish-Chandra homomorphism} maps
$\gamma\colon \D(U/K)\to S(\fa^*)^W$. It can be defined as 
follows. Let $\cU(\fu)$ denote
in the universal enveloping algebra of $\fu$. 
The algebra $\D(U/K)$ is naturally isomorphic
with the quotient $\cU(\fu)^K/\cU(\fu)^K\cap\cU(\fu)\fk$, 
see \cite{GGA}, Thm.\ II.4.6. It follows from 
\cite{GGA} Thm.\ II.5.17 (by application to a symmetric pair
of the non-compact type
with Lie algebras $\fg=\fk+i\fq$ and $\fk$), that there exists
an isomorphism of the quotient 
$\cU(\fu)^\fk/\cU(\fu)^\fk\cap\cU(\fu)\fk$
onto $S(\fa^*)^W$. The Harish-Chandra map results from 
composition of the two,
using that $\cU(\fu)^K\subset\cU(\fu)^\fk$.
We shall need the following fact.

\begin{lemma}\label{l: D(U/K)}
The Harish-Chandra map $\gamma$ is an isomorphism onto $S(\fa^*)^W$.
\end{lemma}

\begin{proof} 
Let $K_0$ denote the identity component of $K$.
It follows from the description of $\gamma$ above,
that it suffices to prove equality between the quotients
$\cU(\fu)^K/\cU(\fu)^K\cap\cU(\fu)\fk$ and
$\cU(\fu)^\fk/\cU(\fu)^\fk\cap\cU(\fu)\fk=
\cU(\fu)^{K_0}/\cU(\fu)^{K_0}\cap\cU(\fu)\fk$.

We shall employ \cite{GGA}, Cor.\ II.4.8, according to which
the two quotients are in bijective linear correspondence with
$S(\fq)^K$ and $S(\fq)^{K_0}$, respectively. It therefore
suffices to prove identity between these two spaces.

Let $p\in S(\fq)^{K_0}$ and let $k\in K$. By means of
the Killing form we regard $p$ as a polynomial function
on $\fq$. The claimed identity amounts to
$p\circ{\rm Ad}k=p$. Notice that $p\circ{\rm Ad}k
\in S(\fq)^{K_0}$,
since $k$ normalizes $K_0$.
According to \cite{GGA}, Cor.\ II.5.12, the elements of
$S(\fq)^{K_0}$ are uniquely determined by restriction to
$\fa$. According to
the lemma below, $k$ is a product of elements from
$K_0$ and $Z_K(\fa)$, and hence it follows that
$p\circ{\rm Ad}k=p$ on $\fa$.
\qed\end{proof}

\begin{lemma}
Each component of $K$ contains an element
from the centralizer $Z_K(\fa)$.
\end{lemma}

\begin{proof} 
Let $k\in K$ be arbitrary. Then ${\rm Ad} k$ maps
$\fa$ to a maximal abelian subspace in $\fq$, hence
to ${\rm Ad} k_0(\fa)$ for some $k_0\in K_0$.
It follows that $k_0^{-1}k$ normalizes $\fa$.
The description of the Weyl group cited in the end of
Section \ref{s: notation} implies that
$N_K(\fa)/Z_K(\fa)=N_{K_0}(\fa)/Z_{K_0}(\fa)$,
hence $k_0^{-1}k\in N_{K_0}(\fa)Z_K(\fa)$
and $k\in K_0Z_K(\fa)$.
\qed\end{proof}

The spherical function
$\psi_\mu$ satisfies the joint eigenequation
\begin{equation}
\label{e: eigenequation D}
D\psi_\mu=\gamma(D,\mu+\rho)\psi_\mu,\qquad D\in\D(U/K)
\end{equation}
(see \cite{BOP} Lemma 2.5). It follows that
$$
(Df)^\sim(\mu)=
\overline{\gamma(D^*,\mu+\rho)}\tilde f(\mu)
$$
where $D^*\in\D(U/K)$ is the adjoint of $D$. 

In particular, the Laplace-Beltrami operator $L$ on $M$
belongs to $\D(U/K)$, and we have
$$
\gamma(L,\lambda)=\langle\gl,\gl\rangle-
\langle\rho,\rho\rangle.
$$
Since $L$ is self-adjoint it follows that
\begin{equation}
\label{e: eigenequation L}
(Lf)^\sim(\mu)=(\langle\mu+\rho,\mu+\rho\rangle-
\langle\rho,\rho\rangle)\tilde f(\mu)
\end{equation}
for all $f\in C^\infty(U/K)^K$.

\eqsection{The estimate of Opdam}
\label{s: Opdam}

In this section we prove part (i) of Theorem \ref{t: PW}.
The proof is based on the following result. 
Let $\bar\Omega$ be the closure of
$$\Omega=\{ X\in\fa \mid \forall\alpha\in\Sigma: 
|\alpha(X)|<\frac \pi2\}.$$

\begin{thm}%
\label{t: estimate}
{\rm [Opdam]} For each $X\in \bar\Omega$ the map
$$\mu\mapsto \psi_\mu(\Exp X), \quad 
\mu\in\Lambda^+(U/K), $$
has an analytic continuation to
$\fa_\C^*$,
denoted $\gl\mapsto \psi_\gl(\Exp X)$, with the following properties.
There exists a constant $C>0$ such that
\begin{equation}
\label{e: Opdam}
|\psi_\gl(\Exp X)|\leq C \,e^{\max_{w\in W}\Re w\gl(X)}
\end{equation}
for all $\gl\in\fa_\C^*$, $X\in\bar \Omega$. 
Furthermore, the map $X\mapsto \psi_\gl(\Exp X)$ is
analytic, and
\begin{equation}
\label{e: W-inv}
\psi_{w(\gl+\rho)-\rho}(\Exp X)=\psi_{\gl}(\Exp X)
\end{equation}
for all $w\in W$.
\end{thm}

\begin{proof} The existence of the analytic continuation
follows from \cite{Opd} Theorem 3.15, by
identification of $\psi_{\mu}(\Exp X)$ with
$G(\mu+\rho,k;X)$, where $G$ is the function appearing there.
Recall that the root system $R$ in
\cite{Opd} is $2\Sigma$. 
For the shift by $\rho$ and (\ref{e: W-inv}), see \cite{BOP}, 
Lemma 2.5. 
It follows from  \cite{Opd} Theorem 6.1 (2) that 
the analytic extension satisfies (\ref{e: Opdam}).
\qed\end{proof} 

\begin{rem} An analytic extension of $\psi_\mu(\Exp X)$
exists for
$X$ in the larger domain $2\Omega$. This was proved
by Faraut (see \cite{BOP}, p.\ 418) and by
Kr\"otz and Stanton (see \cite{KrSt}). However, 
the estimate (\ref{e: Opdam}) has not been
obtained in this generality.\end{rem}

We can now derive Theorem \ref{t: PW} (i).
The following integration formula holds on $M=U/K$
(see \cite{GGA}, p.\ 190), up to normalization of measures:
$$\int_{M} f(x)\,dx
= \int_{K}\int_{A_*} f(ka\cdot o) \delta(a)\,da\,dk$$
where $A_*$ is the torus $\exp\fa$ in $U$ equipped with 
Haar measure, and where $\delta$ is defined by 
$$\delta(\exp H)= \Pi_{\ga\in\gS^+} |\sin i\ga(H)|^{m_\ga}$$
for $H\in\fa$. It follows that
$$
\tilde f(\mu)= \int_{A_*} f(a\cdot o) \psi_\mu(a^{-1}\cdot o)
\delta(a)\,da.$$
Let $R>0$ be sufficiently small, such that
the restriction of $\exp$ to $B_R(0)$
is injective, then if $r<R$ and $f$ 
is $K$-invariant with support inside
$\Exp \bar B_r(0)$, it follows that
\begin{equation}
\label{e: tilde f(mu)}
\tilde f(\mu)= \int_{B_r(0)\cap\fa} f(\Exp H) \psi_\mu(\Exp(-H))
\delta(\exp H)\,dH.
\end{equation}

Assume in addition that $R\leq\pi/(2\|\alpha\|)$
for all $\alpha\in\Sigma$.
Then $B_r(0)\cap\fa\subset\Omega$ for $r<R$, and
it follows from Theorem \ref{t: estimate} that
$\mu\mapsto\tilde f(\mu)$ allows an analytic 
continuation to $\fa^*_\C$, given by the same formula
(\ref{e: tilde f(mu)}), and 
denoted $\tilde f(\gl)$, such that
\begin{equation}
\label{e: est}
|\tilde f(\gl)| \leq 
C \max_{a\in A_*}\{|f(a\cdot o)\delta(a)|\}\, e^{r\|\Re\gl\|}
\end{equation}
where $C$ is a constant depending on $r$, but not on $f$. 
The derivation of the
polynomial decay  of $\tilde f(\gl)$
in (a) of Definition \ref{d: PW space} is then 
easily obtained from  the estimate (\ref{e: est}), when
applied to the function $L^mf$ with a sufficiently high power of $L$,
by means of (\ref{e: eigenequation L}).

The Weyl group transformation property in
part (b) of  Definition \ref{d: PW space} 
follows immediately from (\ref{e: W-inv}).
Hence we can conclude that $\tilde f(\gl)$
belongs to  $\PW_r(\fa)$.

\eqsection{Uniqueness}
\label{s: carlson}

In this section part (iii) of Theorem 
\ref{t: PW} is proved. The proof is based on
the following simple generalization of Carlson's theorem
(see \cite{Boas} p.\ 153).

\begin{lemma}
\label{l: Carlson}
Let $f\colon\C^n\to\C$ be holomorphic. Assume: 

\smallskip
{\rm (i)} There exist a constant $c<\pi$, and
for each $z\in\C^n$ a constant $C$ 
such that
$$|f(z+\zeta e_i)|\leq Ce^{c|\zeta|}$$
for all $\zeta\in\C$, $i=1,\dots,n$. 

{\rm (ii)} $f(k)=0$ for all
$k\in(\Z^{+})^n$.

\smallskip\noindent Then $f=0$.
\end{lemma}

\begin{proof} For $n=1$ this is Carlson's theorem.
In general
it follows by induction that
$z\mapsto f(z,\kappa)$ is identically $0$ on $\C^{n-1}$
for each $\kappa\in\Z^+$.
By a second application of Carlson's theorem it then
follows that $f(z,\zeta)=0$ for all $(z,\zeta)\in\C^n$. 
\qed\end{proof}

It follows that if $X$ is sufficiently close to 0, then
the analytic continuation $\gl\mapsto \psi_\gl(\Exp X)$
in Theorem \ref{t: estimate} is unique, 
when (\ref{e: Opdam}) is required. More precisely,
let $\mu_1,\dots,\mu_n\in\fa^*_\C$ be such that
$\Lambda^+(U/K)=\Z^+\mu_1+\dots+\Z^+\mu_n$.
If $U$ is simply connected the elements
$\mu_1,\dots,\mu_n\in\fa^*_\C$ are the fundamental
weights determined by 
$$\frac{\langle \mu_i,\ga_j\rangle}{\langle\ga_j,\ga_j\rangle}
=\delta_{ij}$$
for the simple roots $\ga_1,\dots,\ga_n$ of $\Sigma^+$.
If $U$ is not simply connected, the $\mu_i$ are suitable 
integral multiples of the fundamental weights,
in order that they correspond to representations of $U$.
If $\|X\|<\pi/\|\mu_i\|$ for all $i$
the uniqueness of the analytic continuation 
now follows by application of
Lemma \ref{l: Carlson}
to the function $z\mapsto f(z_1\mu_1+\dots+z_n\mu_n)$. 

In the same fashion, if $R\leq\pi/\|\mu_i\|$ for all $i$,
it follows from Lemma \ref{l: Carlson}
that for $r<R$
the elements $\varphi\in\PW_r(\fa)$
are uniquely determined on $\Lambda^+(U/K)$, as claimed in Theorem 
\ref{t: PW} (iii).

Notice that the minimal value of $\pi/\|\mu_i\|$ can be 
strictly smaller than the injectivity radius of $U/K$. 
See Remark \ref{r: best R}.

\eqsection{The theorem of Gonzalez}
\label{s: central}

In this section we treat the special case, where the symmetric
space is the compact semisimple Lie group $U$ itself, viewed
as a symmetric space for the product group
$U\times U$ with the action given by $(g,h)\cdot x=gxh^{-1}$.
The stabilizer at $e$ is the diagonal subgroup 
$\Delta=\{(x,x)\mid x\in U\}$ in $U\times U$,
and the corresponding involution of $U\times U$ is 
$(x,y)\mapsto (y,x)$. The $\Delta$-invariant functions
on $U$ are the class functions
(also called central functions), that is, those for which
$f(uxu^{-1})=f(x)$ for all $u,x\in U$. In this case the
local Paley-Wiener theorem was obtained by Gonzalez \cite{Gon}. 
Let us recall his result. 

As before, we denote by $\fh$ a Cartan subalgebra of $\fu$,
and by $\gL^+(\fh)\subset i\fh^*$ the set of dominant
integral linear functionals. For $\mu\in\gL^+(U)$
we denote by $\chi_\mu$ the character of $\pi_\mu$, that is,
$\chi_\mu(x)$ is the trace of $\pi_\mu(x)$ for $x\in U$.
The function $d(\mu)^{-1}\chi_\mu$ is normalized so that
its value at $e$ is 1, and when
$U$ is viewed as a symmetric space, this class function
is exactly the spherical function associated with $\pi_\mu$.
It is however more convenient to use the unnormalized function
$\chi_\mu$ in the definition of the Fourier transform,
since it is a unit vector in $L^2$ (with the
normalized Haar measure on $U$).

Following custom, we thus define the Fourier transform by
$$\hat F(\mu)= \langle F,\chi_\mu \rangle
=\int_U F(u)\overline{\chi_\mu(u)}\,du, \quad \mu\in\gL^+(U),$$
for class functions $F\in L^2(U)^U$. The corresponding
Fourier series is given by
\begin{equation}
\label{e: U Fourier series}
\sum_{\mu\in\gL^+(U)} \hat F(\mu)\chi_\mu(x).
\end{equation}
It converges to $F$ in $L^2$. If $F$ is smooth
it also converges absolutely and uniformly (see \cite{GGA} p.\ 534).

The theorem of Gonzalez \cite{Gon} now reads as follows.
Let $R>0$ be the injectivity radius of $U$.
If $U$ is simply connected, this means that
$R=2\pi/\|\ga\|$ where $\ga$ is the longest root in $\Delta$
(see \cite{Sig} p.\ 318).

\begin{thm}\label{t: Gonzalez}{\rm [Gonzalez]}
Let a class function $F\in C^\infty(U)^U$ be given, 
and let $0<r<R$. Then $F$ belongs to
$C^\infty_r(U)^U$  if and only if the 
Fourier transform $\mu\mapsto \hat F(\mu)$
extends to a holomorphic function $\Phi$
on $\fh^*_\C$ with the following properties

\smallskip
\item{\rm (a)} For each $k\in\N$ there exists a constant $C_k>0$ such that
$$|\Phi(\gl)|\leq C_k(1+\|\gl\|)^{-k} e^{r\|\Re\gl\|}$$
for all $\gl\in\fh_\C^*.$

\item{\rm (b)} $\Phi(w(\gl+\rho)-\rho)=\det(w)\Phi(\gl)$ for
all $w\in W$, $\gl\in\fh_\C^*$.
\end{thm}

Notice that as before the extension $\Phi$ is unique if 
$r$ is sufficiently small. In that case, the Fourier transform,
followed by holomorphic extension, is then a bijection onto the 
space of holomorphic functions satisfying (a) and (b).

\eqsection{Construction of $K$-invariant functions}
\label{s: K-invariant}

The following result is important for the proof of
Theorem \ref{t: PW}.

\begin{lemma}
Let $F\in C^\infty(U)^U$ and define $f\colon U\to\C$
by
$$f(u)=\int_K F(ku) \,dk=\int_K F(uk)\, dk.$$
Then $f\in C^\infty(U/K)^K$ and
\begin{equation}
\label{e: Fourier coefficients}
d(\mu)\tilde f(\mu)=\hat F(\gl)
\end{equation}
for all $\mu\in\gL^+(U/K)$, where $\gl\in\fh^*_\C$
is the extension of $\mu$
determined by $$\gl|_\fa=\mu,\quad\gl|_{\fh_m}=0.$$
\end{lemma}

\begin{proof}
The fact that $f\in C^\infty(U/K)^K$ is clear. From the
uniform convergence of the Fourier series 
(\ref{e: U Fourier series}) it follows
that
$$f(u)=
\sum_{\gl\in\gL^+(U)} \hat F(\gl) \int_K\chi_\gl(uk)\,dk.
$$
By the lemma below we then obtain
$$f(u)=
\sum_{\mu\in\gL^+(U/K)} \hat F(\gl) \psi_\mu(u)
$$
where $\gl$ is the extension of $\mu$ as above.
The statement (\ref{e: Fourier coefficients}) now follows
by comparison with (\ref{e: U/K Fourier series}).\qed
\end{proof}

\begin{lemma}
Let $\gl\in\gL^+(U)$ and $\mu=\gl|_\fa$.
If $\gl\in\gL^+_K(U)$ then
$$\int_K\chi_\gl(uk)\,dk=\psi_\mu(u)$$ for
all $u\in U$, and otherwise 
$
\int_K\chi_\gl(uk)\,dk=0$.
\end{lemma}

\begin{proof} (See also \cite{GGA}, p.\ 417).
The function $u\mapsto\int_K\chi_\gl(uk)\,dk$
is a $K$-fixed vector in the right representation 
generated by $\chi_\gl$, which is equivalent
with $\pi_\gl$,
hence it vanishes
if $\gl\notin \gL^+_K(U)$.

Assume $\gl\in \gL^+_K(U)$, and choose an orthonormal basis
$v_1,\dots, v_d$ for the representation space $V$,
such that $v_1$ is $K$-fixed. Then 
$$\int_K \chi_\gl(uk)\,dk= 
\int_K
\sum_{i=1}^{d}\langle\pi_\gl(u)\pi_\gl(k)v_i,v_i\rangle \,dk
. $$
Since the
operator $\int_K\pi_\gl(k)\,dk$ is the orthogonal projection onto
$V^K$, it follows that
$$
\int_K \chi_\gl(uk)\,dk=
\langle\pi_\gl(u)v_1,v_1\rangle=
\psi_\mu(u). 
$$
as claimed.\qed
\end{proof}

\begin{lemma} 
\label{l: support}
Let $F\in C^\infty(U)^U$ and $f\colon U/K\to\C$
be as above.
If $F\in C_r^\infty(U)^U$ for some $r>0$
then $f\in C_r^\infty(U/K)^K$. 
\end{lemma}

\begin{proof} 
Let $x\in M$ with $f(x)\neq 0$ and choose $X\in\fq$ such that
the curve on $M$ given by
$t\mapsto\gamma(t)=\Exp(tX)$ where
$t\in[0,1]$, is a
minimal geodesic from $o$ to $x$. 
The length of $\gamma$ is $\|X\|$.
 
Let $x=u\cdot o$ where $u\in U$, then there exists 
$k\in K$ such that $F(ku)\neq 0$.
Hence $ku=\exp Y$ where $Y\in\fu$ with $\|Y\|<r$.
Let $Z=\Ad(k^{-1})Y$, then $\|Z\|=\|Y\|<r$.
The smooth curve $\xi(t)=\exp(tZ)\cdot o,$
where $t\in[0,1]$, also joins $o$ to $x$. 
Hence it has length $\ell(\xi)\geq\|X\|$.

Let $L_u$ denote left translation by $u$,
then $\xi'(t)=dL_{\exp(tY)}(\xi'(0))$, and hence 
$\|\xi'(t)\|=\|\xi'(0)\|$ for all $t$. Let
$Z_\fq$ denote the $\fq$ component
of $Z$ in the orthogonal decomposition $\fu=\fk+\fq$.
Then $\xi'(0)=Z_\fq$, and we conclude that 
$$\|X\|\leq \ell(\xi)=
\int_0^1 \|\xi'(t)\|\,dt=\|Z_\fq\|\leq \|Z\|=\|Y\|<r.$$
Thus $f$ has support in $\Exp\bar B_r(0)$\qed
\end{proof}

\eqsection{The result of Rais}
\label{s: Rais}

The following result is due to M. Rais.  
Let $r>0$ and recall that a
holomorphic function $\varphi$ on  $\fh_\C^*$ is said
to be of exponential type $r$ if it satisfies (a) of
Theorem \ref{t: Gonzalez}. 
Let $\tilde W$ denote the Weyl group of the 
root system $\Delta$ on $\fh$.
Let $l=|\tilde W|$, and let
$P_1,\dots,P_l$ be a basis for $S(\fh^*)$ over
$I(\fh^*)=S(\fh^*)^{\tilde W}$ (see \cite{GGA}, p.\ 360).

\begin{thm} 
\label{t: Rais}
For each holomorphic function $\psi$ 
of exponential type $r$
there exist unique $\tilde W$-invariant 
holomorphic functions $\phi_1,\dots,\phi_l$
of exponential type $r$
such that $\psi=P_1\phi_1+\dots+P_l\phi_l$.
\end{thm}

\begin{proof}
See \cite{CD2}, Appendix B.\qed
\end{proof}

In the following statement, we regard $\fa^*$ as a
subset of $\fh^*$, by trivial extension on $\fh_m$.
Likewise $\fh_m^*$ is regarded as a subspace
by trivial extension on $\fa$. Then 
$\fh^*=\fa^*\oplus\fh_m^*$ holds as an orthogonal
sum decomposition.

\begin{Cor} 
\label{c: Rais}
There exist a collection of polynomials
$p_1\dots,p_l\in S(\fa^*)^W$ with the following property.
For each $W$-invariant 
holomorphic function $\varphi$ on $\fa^*_\C$ 
of exponential type $r$,
there exist $\tilde W$-invariant 
holomorphic functions $\phi_1,\dots,\phi_l$ on $\fh^*_\C$
of exponential type $r$,
such that 
\begin{equation}
\label{e: Rais}
\varphi=p_1(\phi_1|_{\fa^*_\C})+\dots+p_l(\phi_l|_{\fa^*_\C}).
\end{equation}
\end{Cor}

\begin{proof}(See also \cite{FJ} p.\ 30). Notice that
when $\phi_j$ is $\tilde W$-invariant, then
$\phi_j|_{\fa^*}$ is
$W$-invariant, since the normalizer
in $\tilde W$ of $\fa$ maps surjectively onto $W$
(see \cite{GGA} p. 366).

Fix a holomorphic function $\varphi_m$ 
on $\fh_{m\C}^*$ of exponential type
$r$, with the value $\varphi_m(0)=1$.
Put $\psi(\gl)=\varphi(\gl_1)\varphi_m(\gl_2)$, where
$\gl_1$ and $\gl_2$ are the components of $\gl$.
Then $\psi$ is of exponential type $r$, and we can apply
Theorem \ref{t: Rais}. 
The restriction of $\psi$
to $\fa^*$ is exactly $\varphi$.
Taking restrictions to $\fa^*$ we thus obtain 
(\ref{e: Rais})  with $p_j=P_j|_{\fa^*}$.
The desired expression is obtained by averaging over~$W$.
\qed\end{proof}

\eqsection{Proof of the main theorem}
\label{s: surjective}

It remains to be seen that every function 
$\varphi\in \PW_r(\fa)$ is the extension of $\tilde f$
for some $f\in C^\infty(U/K)^K$. 

Thus let $\varphi\in \PW_r(\fa)$ be given. Let
$p_1,\dots,p_l$ and $\phi_1,\dots,\phi_l$ be as in 
Corollary \ref{c: Rais}, applied to the 
$W$-invariant function
$\gl\mapsto \varphi(\gl-\rho)$ on $\fa^*$.
By Lemma \ref{l: D(U/K)}
there exist $D_j\in \D(U/K)$ such that
$\overline{\gamma(D_j^*,\gl)}=p_j(\gl)$
for $\gl\in i\fa^*$. 

It follows from the Weyl dimension formula
(\cite{GGA} p.\ 502)
that $\mu\mapsto d(\mu)$ extends
to a polynomial on $\fh^*$ which satisfies
the transformation property (b) of Theorem
\ref{t: Gonzalez}.
Hence the function on $\fh^*_\C$ defined by
$\Phi_j(\gl)=d(\gl)\phi_j(\gl+\rho)$
satisfies both (a) and (b) in that theorem,
and thus we can find $F_1,\dots,F_n\in C^\infty_r(U)^U$
such that $$\hat F_j(\mu)=\Phi_j(\mu)$$ for all $\mu$.

Let $f_j(uK)=\int_K F_j(uk)\, dk$
and define $f=\sum D_j f_j$. 
Then by Lemma \ref{l: support}
we have $f\in C^\infty_r(U/K)^K$, and it follows from
(\ref{e: Fourier coefficients}) that
\begin{eqnarray*}
\tilde f(\mu)&=&\sum_j \gamma(D_j,\mu+\rho)\tilde f_j(\mu)\\
&=&\sum_j \gamma(D_j,\mu+\rho)d(\mu)^{-1}\hat F_j(\mu)\\
&=&\sum_j p_j(\mu+\rho)\phi_j(\mu+\rho)=\varphi(\mu).\qquad\square
\end{eqnarray*}

\eqsection{The sphere $S^2$}
\label{s: sphere}

Let $M=S^2=\{(x,y,z)\in\real^3\mid x^2+y^2+z^2=1\}$,
then $M$ can be realized as a homogeneous space
for $U=\SU(2)$ with the following action.
Identify $\real^3$ with the space of Hermitian
$2\times 2$ matrices $H$ with trace 0,
$$H=\pmatrix{z&x+iy\cr x-iy&-z\cr},$$
then $u.H=uHu^{-1}$ for $u\in U$. The
stabilizer of the point $o=(0,0,1)\in S^2$ is 
the set of diagonal elements in $U$, and the
diagonal element
$$k_\theta=\pmatrix{e^{i\theta/2}&0\cr 0&e^{-i\theta/2}\cr}$$
acts by rotation around the $z$-axis
of angle $\theta$.

A $K$-invariant function on $M$ is determined
by its values along the elements $(x,y,z)=(0,\sin t,\cos t)$,
and it thus becomes
identified as an even function of $t\in [-\pi,\pi]$.
With the notation of above, the 
function is identified through the map $t\mapsto f(a_t\cdot o)$ where
$$a_t=
\pmatrix{\cos(t/2) & i\sin(t/2)\cr 
i\sin(t/2)&\cos(t/2)\cr}.$$

The irreducible representations of $U$ are parametrized
by half integers $l=0,\frac12,1,\frac32,\dots$, 
where $\pi_l$ has dimension $2l+1$, and the spherical
representations are those for which $l$ is an integer.
The corresponding spherical functions are given by
$\psi_l(a_t)=P_l(\cos t)$, where $P_l$ is the
$l$'th Legendre polynomial.

The Fourier series of a $K$-invariant function on $S^2$,
identified as an even function on $[-\pi,\pi]$,
is then the Fourier-Legendre series
$$\sum_{l=0}^\infty (2l+1) \tilde f(l) P_l(\cos t)$$
where 
$$\tilde f(l)= \frac12
\int_0^\pi f(t) P_l(\cos t) \sin t\,dt.$$

Our local Paley-Wiener theorem asserts
the following for $r<\pi$: 

\smallskip
{\it An even function
$f\in C^\infty(-\pi,\pi)$ is supported in $[-r,r]$
if and only if the Legendre transform
$l\mapsto \tilde f(l)$ of $f$ extends to
an entire function $g$ on $\C$
of exponential type
$$|g(\gl)|\leq C_k(1+|\gl|)^{-k}e^{r|\Im\gl|}$$
such that $g(\gl-\frac12)$ is an even function of $\gl$.
The extension $g$ with these properties is unique.

Moreover, every such function $g$ on $\C$ is obtained
in this fashion from a unique function
$f\in C_r^\infty(-\pi,\pi)$.}

\smallskip
Essentially this is the result stated
by Koornwinder (and attributed to Beurling) in
\cite{TomK} p.\ 158.



\def\adritem#1{\hbox{\small #1}}
\def\distance{\hbox{\hspace{3.5cm}}}
\def\addgestur{\vbox{
\adritem{G.~\'Olafsson}
\adritem{Department of Mathematics}
\adritem{Louisiana State University}
\adritem{Baton Rouge}
\adritem{LA 70803}
\adritem{USA}
\adritem{E-mail: olafsson@math.lsu.edu}
}
}
\def\addhenrik{\vbox{
\adritem{H.~Schlichtkrull}
\adritem{Department of Mathematical Sciences}
\adritem{University of Copenhagen}
\adritem{Universitetsparken 5}
\adritem{2100 K\o benhavn \O}
\adritem{Denmark}
\adritem{E-mail: schlicht@math.ku.dk}
}
}
\mbox{\ }
\vfill
\hbox{\vbox{\addgestur}\vbox{\distance}\vbox{\addhenrik}}
\end{document}